\documentclass[a4paper,12pt]{article}
\usepackage{amssymb}
\usepackage[pagewise, mathlines, displaymath, running]{lineno}

\newcommand{\RR}{ I\!\!R}

\newcommand{\R}{\mathbb{R}}

\newcommand{\N}{\mathbb{N}}

\newcommand{\lambd}{\mu}

\newcommand{\beq}{\begin{equation} }
\newcommand{\eqq}{\end{equation} }
\newcommand{\cuad}{{\sqcap\kern-.68em\sqcup}}

\newcommand{\equ}[1]{(\ref{#1})}

\newtheorem{teo}{Theorem}[section]

\newtheorem{proposition}{Proposition}[section]

\newtheorem{lemma}{Lemma}[section]

\newtheorem{remark}{Remark}[section]
\newcommand{\bremark}{\begin{remark} \em}
\newcommand{\eremark}{\end{remark} }

\def\beeq{\begin{equation}}
\def\eeq{\end{equation}}
\newcommand{\begeqaet}{\begin{eqnarray*}}
\newcommand{\eneqaet}{\end{eqnarray*}}

\usepackage{graphicx}

\begin{document}

\begin{center}{\bf  \large Self-generated interior blow-up solutions in fractional elliptic equation with absorption}\medskip

\bigskip
\bigskip

{Huyuan Chen,  \  Patricio Felmer \\~\\}

 Departamento de Ingenier\'{\i}a  Matem\'atica and
Centro de Modelamiento Matem\'atico
 UMR2071 CNRS-UChile,
 Universidad de Chile\\
 Casilla 170 Correo 3, Santiago, Chile.\\
 {\sl  (hchen@dim.uchile.cl, \ pfelmer@dim.uchile.cl)}
 \bigskip
 \\and\\
 \bigskip
 { Alexander Quaas \\~\\}
Departamento de Matem\'{a}tica, Universidad T\'{e}cnica Federico
Santa Mar\'{i}a\\ Casilla: V-110, Avda. Espa\~{n}a 1680,
Valpara\'{i}so, Chile\\
 {\sl  (alexander.quaas@usm.cl)}

\bigskip

\bigskip
\begin{abstract}
In this paper we study positive solutions to problem involving the
fractional Laplacian
\begin{equation}\label{0000}
\left\{ \arraycolsep=1pt
\begin{array}{lll}
 (-\Delta)^{\alpha} u(x)+|u|^{p-1}u(x)=0,\ \ \ \ &
x\in\Omega\setminus\mathcal{C},\\[2mm]
u(x)=0,\ &x\in\Omega^c,\\[2mm]
\lim_{x\in\Omega\setminus\mathcal{C}, \
x\to\mathcal{C}}u(x)=+\infty,
\end{array}
\right.
\end{equation}
where $p>1$ and $\Omega$ is an open bounded $C^2$ domain in
$\mathbb{R}^N$,  $\mathcal{C}\subset \Omega$ is a compact $C^2$
manifold with $N-1$ multiples dimensions and without boundary, the
operator $(-\Delta)^{\alpha}$ with $\alpha\in(0,1)$ is the
fractional Laplacian.

We consider the existence of positive solutions for problem
(\ref{0000}). Moreover, we further analyze uniqueness, asymptotic
behaviour and nonexistence.

\end{abstract}

\vspace{1mm}
\end{center}
 \setcounter{equation}{0}
\section{Introduction}

In 1957, a fundamental contribution due to Keller in \cite{K} and
Osserman in \cite{O} is the study of boundary  blow-up solutions for the non-linear elliptic equation
\begin{equation}\label{1.1.2}
\left\{ \arraycolsep=1pt
\begin{array}{lll}
 -\Delta u+h(u)=0\  \
 \rm{in}\  \ \Omega,\\[2mm]
\lim_{x\in \Omega, x\to\partial\Omega}u(x)=+\infty.    
\end{array}
\right.
\end{equation}
They proved the existence of solutions to
(\ref{1.1.2}) when $h:\R\to[0,+\infty)$ is a locally Lipschitz
continuous function which is nondecreasing and satisfies the so called
Keller-Osserman condition. From then on, the result of Keller and Osserman  has been extended by numerous mathematicians in
various ways, weakening the assumptions on the domain, generalizing the differential operator and the
nonlinear term for equations and systems. The case of  $h(u)=u_+^p$ with $p=\frac{N+2}{N-2}$
 is studied by Loewner and Nirenberg \cite{LN5}, where in particular uniqueness and asymptotic behavior were obtained. After that,  Bandle and Marcus \cite{BM1}  obtained uniqueness and asymptotic for more general non-linearties $h$. Later, Le Gall in
 \cite{G5} established a uniqueness result of problem (\ref{1.1.2})
in the  domain whose boundary is non-smooth when  $h(u)=u_+^2$.
Marcus and V\'{e}ron \cite{MV2,MV1} extended the uniqueness of blow-up
solution for (\ref{1.1.2}) in general domains whose boundary is
locally represented as a graph of a continuous function when
$h(u)=u_+^p$ for $p>1$. Under this special assumption on $h$, Kim
\cite{K1} studied the existence and uniqueness of  boundary blow-up
solution to (\ref{1.1.2}) in bounded domains $\Omega$ satisfying
$\partial \Omega=\partial \bar\Omega$. For another interesting contributions to
boundary blow-up solutions see for example Kondratev, Nikishkin \cite{KN},   Lazer, McKenna \cite{LM}, Arrieta and Rodr\'{\i}guez-Bernal \cite{AR}, Chuaqui, Cort\'{a}zar,  Elgueta and J. Garc\'{i}a-Meli\'{a}n \cite{CCE},
del Pino and Letelier \cite{DL},  D\'{i}az and  Letelier \cite{DL1}, Du and  Huang \cite{DH},  Garc\'{i}a-Meli\'{a}n  \cite{G}, V\'eron \cite{V}, and the reference therein.

In a recent work, Felmer and Quaas \cite{FQ0} considered a version of Keller and Osserman problem for a class of non-local operator.
Being more precise, they considered as a particular case the fractional elliptic problem

\begin{equation}\label{0.1.1}
\left\{ \arraycolsep=1pt
\begin{array}{lll}
 (-\Delta)^{\alpha} u(x)+|u|^{p-1}u(x)=f(x),\ \ \ \ &
x\in\Omega,\\[2mm]
u(x)=g(x),& x\in\bar\Omega^c,\\[2mm]
\lim_{x\in\Omega,\ x\to\partial\Omega}u(x)=+\infty,
\end{array}
\right.
\end{equation}
where $p>1$, $f$ and $g$ are appropriate functions and $\Omega$ is a bounded domain with $C^2$ boundary.
The operator $(-\Delta)^{\alpha}$ is the fractional Laplacian which
is defined as
\begin{equation}\label{fractional}
(-\Delta)^\alpha u(x)=-\frac12
\int_{\R^N}\frac{\delta(u,x,y)}{|y|^{N+2\alpha}}dy,\ \ x\in \Omega,
\end{equation}
with  $\alpha\in(0,1)$ and $\delta(u,x,y)=u(x+y)+u(x-y)-2u(x)$.

In \cite{FQ0} the authors proved the existence of a solution  to
(\ref{0.1.1}) provided that $g$ explodes at the boundary and
satisfies other technical conditions. In case the function $g$ blows
up with an explosion rate as $d(x)^{\beta}$, with $\beta\in
[-\frac{2\alpha}{p-1},0)$ and $d(x)=dist(x,\partial\Omega)$, it is shown that the
solution satisfies
$$
0< \liminf_{x\in\Omega, x\to\partial\Omega}u(x)d(x)^{-\beta}\le
\limsup_{x\in\Omega,
x\to\partial\Omega}u(x)d(x)^{\frac{2\alpha}{p-1}}<+\infty.
$$
Here
 the explosion is driven by the external value  $g$ and
the external source $f$ has a secondary role, not
intervening in the explosive character of the solution.

More recently,
Chen, Felmer and Quaas \cite{CFQ} extended the results in \cite{FQ0} studying existence, uniqueness and non-existence of  boundary blow-up
solutions when the function $g$ vanishes and the explosion on the boundary is driven by the external source $f$, with weak or strong  explosion rate.  Moreover,  the results are extended even to the case where the boundary blow-up solutions in driven internally, when the external source and value, $f$ and $g$, vanish. Existence, uniqueness, asymptotic behavior  and non-existence results for blow-up solutions of \equ{0.1.1} are considered in \cite{CFQ}. In the analysis developed in \cite{CFQ}, a key role is played by the function $C:(-1,0]\to\R,$ that governs the behavior of the solution near the boundary. The function $C$ is defined as
\begin{equation}\label{1.4}
C(\tau)=\int^{+\infty}_{0}\frac{\chi_{(0,1)}(t)|1-t|^{\tau}+(1+t)^{\tau}-2}{t^{1+2\alpha}}dt
\end{equation}
and it possess exactly one zero in $(-1,0)$ and we call it $\tau_0(\alpha)$. In what follows we explain with more details the results in the case of vanishing external source and values, that is $f=0$ in $\Omega$ and $g=0$ in
$\bar\Omega^c$, which is the  case we will consider in this paper. In Theorem 1.1 in  \cite{CFQ}, we proved that
whenever
$$
1+2\alpha<p<1-\frac{2\alpha}{\tau_0(\alpha)},
$$
then problem
 (\ref{0.1.1}) admits a unique positive solution $u$ such that
$$
0<\liminf_{x\in\Omega,x\to\partial\Omega}u(x)d(x)^{\frac{2\alpha}{p-1}}
\leq\limsup_{x\in\Omega,x\to\partial\Omega}u(x)d(x)^{\frac{2\alpha}{p-1}}<+\infty.
$$
On the other hand, we proved that when
 $p\ge1,$ then problem
 (\ref{0.1.1}) does not admit any solution $u$ such that
\begin{equation}\label{1.5}
0<\liminf_{x\in\Omega,x\to\partial\Omega}u(x)d(x)^{-\tau}
\leq\limsup_{x\in\Omega,x\to\partial\Omega}u(x)d(x)^{-\tau}<+\infty,
\end{equation}
for any $\tau\in(-1,0)\setminus\{\tau_0(\alpha),-\frac{2\alpha}{p-1}\}$. We observe that the non-existence result does not include the  case when $u$ has an asymptotic behavior of the form $d(x)^{\tau_0(\alpha)}$, where $\tau_0(\alpha)$ is precisely where $C$ vanishes. We have a a special existence result in this case, precisely
if
$$
\max\{1-\frac{2\alpha}{\tau_0(\alpha)}+\frac{\tau_0(\alpha)+1}{\tau_0(\alpha)},1\}<p<1-\frac{2\alpha}{\tau_0(\alpha)},
$$
then, for any $t>0$, problem
 (\ref{0.1.1}) admits   a positive solution $u$  such that
$$
\lim_{x\in\Omega,x\to\partial\Omega}u(x)d(x)^{-\tau_0(\alpha)}=t.
$$

 Motivated by
these results and in view of the  non-local character of the fractional Laplacian we are interested in another class of blow-up solutions. We want to study solutions that vanish at the boundary of the domain $\Omega$ but that  explodes at the interior of the domain, near a prescribed embedded manifold. From now on, we assume that $\Omega$ is an open bounded domain in $\R^N$ with $C^2$ boundary, and that there is a $C^2$, $(N-1)$-dimensional manifold $\mathcal{C}$ without boundary, embedded in $\Omega$, such that,  it separates $\Omega\setminus \mathcal{C}$ in exactly two connected components. We denote by
 $\Omega_1$ the inner component and  by $\Omega_2$ the  external component, that is $\bar\Omega_1\cap \partial\Omega=\emptyset$ and $\bar\Omega_2\cap \partial\Omega=\partial\Omega.$  Throughout the paper we will consider the distance function
 \begin{equation}\label{distance}
 D:\Omega\setminus \mathcal{C}\to \R_+, \quad D(x)={\rm dist}(x, \mathcal{C}).
  \end{equation}

Let us consider the equations, for $i=1,2$,
 \begin{equation}\label{ee1.1}
\left\{ \arraycolsep=1pt
\begin{array}{lll}
 (-\Delta)^{\alpha} u(x)+|u|^{p-1}u(x)=0,\ \ \ \ &
x\in\Omega_i,\\[2mm]
u(x)=0,& x\in\bar\Omega_i^c,\\[2mm]
\lim_{x\in\Omega_i,\ x\to\partial\Omega_i}u(x)=+\infty,
\end{array}
\right.
\end{equation}
which have solutions $u_1$ and $u_2$, for $i=1, 2$ respectively, in the appropriate range of the parameters.
In the local case, that is, $\alpha=1$, these two solutions certainly do not interact among each other, but when $\alpha\in (0,1)$, due to the non-local character of the fractional Laplacian and the non-linear character of the equation the solutions on each side of $\Omega$ interact and it is precisely the purpose of this paper to study their existence, uniqueness and non-existence.

In precise terms we consider the equation
\begin{equation}\label{inside}
\left\{ \arraycolsep=1pt
\begin{array}{lll}
 (-\Delta)^{\alpha} u(x)+|u|^{p-1}u(x)=0,\ \ \ \ &
x\in\Omega\setminus\mathcal{C},\\[2mm]
u(x)=0,\ &x\in\Omega^c,\\[2mm]
\lim_{x\in\Omega\setminus\mathcal{C}, \
x\to\mathcal{C}}u(x)=+\infty,
\end{array}
\right.
\end{equation}
where $p>1$,  $\Omega$ and  $\mathcal{C}\subset \Omega$ are as described above. The explosion of the solution near $\mathcal{C}$ is governed by a function $c:(-1,0]\to \R,$
 defined as
 \begin{equation}\label{1.2}
c(\tau)=\int_0^{+\infty}\frac{|1-t|^\tau+(1+t)^\tau-2}{t^{1+2\alpha}}dt.
\end{equation}
This function plays the role of the function $C$ used in \cite{CFQ}, but it has certain differences. In Section \S 2 we  prove the existence of a number  $\alpha_0\in (0,1)$ such that  $\alpha\in [\alpha_0,1)$ the function $c$ is always positive in $(-1,0)$, while if $\alpha\in (0, \alpha_0)$ then there exists  exists a  unique
$\tau_1(\alpha)\in(-1,0)$ such that $c(\tau_1(\alpha))=0$ and $c(\tau)>0$ in $(-1,\tau_1(\alpha))$ and $c(\tau)<0$ in $(\tau_1(\alpha), 0)$, see Proposition \ref{lm 2.1}.
We notice here that $\tau_1(\alpha)>\tau_0(\alpha)$ if $\alpha\in (0, \alpha_0)$.

Now we are ready to state our main theorems on the existence uniqueness and
asymptotic behavior of interior  blow-up solutions to equation (\ref{inside}). These theorems deal separately the case $\alpha\in (0,\alpha_0)$ and $\alpha\in [\alpha_0,1)$.
\begin{teo}\label{th 1.2}
Assume that $\alpha\in (0, \alpha_0)$ and the  assumptions on $\Omega$ and $\mathcal{C}$. Then we have: \\
$(i)$ If
\begin{equation}\label{0.1.21}
1+2\alpha<p<1-\frac{2\alpha}{\tau_1(\alpha)},
\end{equation} then
problem (\ref{inside}) admits a unique positive solution $u$
satisfying
\begin{equation}\label{0.1.2}
0<\liminf_{x\in\Omega\setminus\mathcal{C},x\to\mathcal{C}}u(x)D(x)^{\frac{2\alpha}{p-1}}
\leq\limsup_{x\in\Omega\setminus\mathcal{C},x\to\mathcal{C}}u(x)D(x)^{\frac{2\alpha}{p-1}}<+\infty.
\end{equation}

\noindent$(ii)$ If
\begin{equation}\label{observartion}
\max\{1-\frac{2\alpha}{\tau_1(\alpha)}+\frac{\tau_1(\alpha)+1}{\tau_1(\alpha)},1\}<p<1-\frac{2\alpha}{\tau_1(\alpha)}.
\end{equation}
Then, for any $t>0$, there is  a  positive solution $u$ of problem (\ref{inside})
satisfying
\begin{equation}\label{o3.1.5}
\lim_{x\in \Omega\setminus{\mathcal C}, x\to{\mathcal C}}u(x)D(x)^{-\tau_1(\alpha)}=t.
\end{equation}
$(iii)$ If one of the following three conditions holds
\begin{itemize}
\item[\rm{a)}] $1<p\le 1+2\alpha$ and $\tau\in(-1,0)\setminus\{ \tau_1(\alpha) \}$,
\item[\rm{b)}] $1+2\alpha<p< 1-\frac{2\alpha}{\tau_1(\alpha)}$  and $\tau\in(-1,0)\setminus\{\tau_1(\alpha),-\frac{2\alpha}{p-1}\}$ or
\item[\rm{c)}] $p\ge 1-\frac{2\alpha}{\tau_1(\alpha)}$ and $\tau\in (-1,0)$,
\end{itemize}
then problem
(\ref{inside}) does not admit any  solution $u$ satisfying
\begin{equation}\label{0.1.3}
0<\liminf_{x\in\Omega\setminus\mathcal{C},x\to\mathcal{C}}u(x)D(x)^{-\tau}
\le\limsup_{x\in\Omega\setminus\mathcal{C},x\to\mathcal{C}}u(x)D(x)^{-\tau}<+\infty.
\end{equation}
\end{teo}
We observe that this theorem is similar to Theorem 1.1 in \cite{CFQ}, where the role of $\tau_0(\alpha)$ is played here by $\tau_1(\alpha)$. A quite
 different situation occurs when
 $\alpha\in [\alpha_0,1)$ and the function $c$ never vanishes in $(-1,0)$. Precisely, we have
\begin{teo}\label{th 1.1}
Assume that $\alpha\in [\alpha_0,1)$   and the  assumptions on $\Omega$ and $\mathcal{C}$.  Then we have:\\
$(i)$ If $ p>1+2\alpha$,
 then
problem (\ref{inside}) admits a unique positive solution $u$ satisfying  (\ref{0.1.2}). \\
$(ii)$ If $p>1,$
 then problem
(\ref{inside}) does not admit any  solution $u$ satisfying \equ{0.1.3}
for any $\tau\in(-1,0)\setminus\{-\frac{2\alpha}{p-1}\}$.
\end{teo}

Comparing Theorem \ref{th 1.2} with Theorem \ref{th 1.1} we see that the range of existence for the absorption term is quite larger for the second one, no constraint from above.  The main difference with the results in \cite{CFQ}, Theorem 1.1, with vanishing $f$ and $g$ occurs when $\alpha$ is large and the function $c$ does not vanish, allowing thus for existence for all $p$ large.  This difference  comes from the fact that
the fractional Laplacian is a non-local operator so that in the interior blow-up, in each of the domains $\Omega_1$ and $\Omega_2$ there is a non-zero external value, the solutions itself acting on the other side of $\mathcal{C}$.

The proof of our theorems is obtained through the use of super and sub-solutions as in \cite{CFQ}. The main difficulty here is  to find the appropriate super
and sub-solutions to apply the iteration technique to fractional
elliptic problem (\ref{inside}). Here we make use of some precise
estimates based on the function $c$ and the distance function $D$ near $\mathcal{C}$.

 This article is organized as
follows. In section \S2, we introduce some preliminaries and we prove the main estimates of the behavior of the fractional Laplacian when applied to suitable powers of the function $D$.  In section \S3 we prove the existence of
solution to problem (\ref{inside}) as given in
Theorem \ref{th 1.2} and Theorem \ref{th 1.1}. Finally, in Section \S4 we  prove the
uniqueness and nonexistence statements of these theorems.

 \setcounter{equation}{0}
\section{Preliminaries}

In this section, we recall some basic results from \cite{CFQ} and
obtain some useful estimate, which will be used in constructing
super and sub-solutions of problem (\ref{inside}). The first result
states as:
\begin{teo}\label{th 2.1}
Assume that $p>1$ and there are super-solution $\bar U$ and
sub-solution $\underline{U}$ of problem (\ref{inside}) such that
$$
\bar U\geq \underline{U}\ \ {\rm{in}}\ \Omega\setminus\mathcal{C},\quad  \liminf_{x\in\Omega\setminus\mathcal{C},x\to\mathcal{C}}\underline U(x)=+\infty,\quad \bar
U=\underline{U}=0\ \ {\rm{in}}\ \Omega^c.$$ Then problem
(\ref{inside}) admits at least one positive solution $u$ such that
$$\underline{U}\leq u\leq \bar U\ \ {\rm{in}}\
\Omega\setminus\mathcal{C}.$$
\end{teo}
{\bf Proof.}  The procedure is similar to the proof of Theorem 2.6 in \cite{CFQ}, here we give the main differences.

Let us define $\Omega_n:=\{ x\in \Omega \,|\,D(x)>1/n\}$ then we solve
\begin{equation}\label{ecnn}
\left\{ \arraycolsep=1pt
\begin{array}{lll}
 (-\Delta)^{\alpha} u_n(x)+|u_n|^{p-1}u_n(x)=0,\ \ \ \ &
x\in\Omega_n,\\[2mm]
u_n(x)=\underline U,\ &x\in\Omega_n^c.\\[2mm]
\end{array}
\right.
\end{equation}
To find these solutions of \equ{ecnn} we observe that for fix $n$ the method of section 3 of \cite{FQ0} applies even if the domain is not connected since
the estimate of Lemma 3.2 holds with $\delta<1/2n$ (see also Proposition 3.2 part  ii) in  \cite{CFQ}), form here sub and super-solution can be construct for the Dirichlet problem and then existence holds for \equ{ecnn} by an iteration technique (see also section 2 of  \cite{CFQ} for that procedure).
Then as in  Theorem 2.6 in \cite{CFQ} we have $$\underline U\leq u_{n}\leq
u_{n+1}\leq \bar U\ \ \rm{in}\ \ \Omega.$$

By monotonicity of $u_n$, we can define
$$
u(x):=\lim_{n\to+\infty}u_n(x),\ x\in \Omega\ \ {\rm{and}}\ \
u(x):=0,\ x\in \Omega^c.
$$
Which, by a stability property, is a solution of  problem (\ref{inside}) with the desired properties.
\hfill$\Box$\\


In order to prove our existence result, it is crucial to have
available super and sub-solutions to problem (\ref{inside}). To this
end,  we start describing the properties of $c(\tau)$ defined in
(\ref{1.2}), which is a $C^2$ function in $(-1,0)$.
\begin{proposition}\label{lm 2.1}
 There exists a unique $\alpha_0\in(0,1)$ such that\\
$(i)$ For $\alpha\in[\alpha_0,1)$, we have $c(\tau)>0,$ for all $\tau\in(-1,0);$\\[1mm]
$(ii)$ For any $\alpha\in(0,\alpha_0)$, there exists unique
$\tau_1(\alpha)\in(-1,0)$ satisfying
\begin{equation}\label{00}
c(\tau)\ \left\{ \arraycolsep=1pt
\begin{array}{lll}
 >0,\quad &\rm{if}\quad \tau\in(-1,\tau_1(\alpha)),\\[2mm]
 =0,\ &\rm{if}\quad   \tau=\tau_1(\alpha),\\[2mm]
 <0,\ \ \ \ &\rm{if}\quad \tau\in(\tau_1(\alpha),0)
\end{array}
\right.
\end{equation}
and
 \begin{equation}\label{1.3.2}
\lim_{\alpha\to\alpha_0^-}\tau_1(\alpha)=0 \quad\rm{ and}\quad
\lim_{\alpha\to0^+}\tau_1(\alpha)=-1.
\end{equation}
Moreover,
$\tau_1(\alpha)>\tau_0(\alpha)$, for all $\alpha\in(0,\alpha_0)$,
where $\tau_0(\alpha)\in(-1,0)$ is the unique   zero of $C(\tau)$, defined
in (\ref{1.4}).
\end{proposition}
{\bf Proof.} From \equ{1.2}, differentiating twice we find that
\begin{equation}\label{2.01}
c''(\tau)=\int_0^{+\infty}\frac{|1-t|^{\tau}(\log|1-t|)^2+(1+t)^{\tau}(\log(1+t))^2}{t^{1+2\alpha}}dt>0,
\end{equation}
so that $c$ is strictly convex in $(-1,0)$. We also see easily that
\begin{equation}\label{c0}
c(0)=0\quad\mbox{and}\quad \lim_{\tau\to -1^+} c(\tau)=\infty.
\end{equation}

Thus,  if $c'(0)\le0$ then
$c(\tau)>0$ for $\tau\in(-1,0)$ and if $c'(0)>0$, then there exists
$\tau_1(\alpha)\in(-1,0)$ such that $c(\tau)>0$ for
$\tau\in(-1,\tau_1(\alpha))$, $c(\tau)<0$ for
$\tau\in(\tau_1(\alpha),0)$ and $c(\tau_1(\alpha))=0$. In order to complete our proof, we have to analyze the sign of $c'(0)$, which depends on $\alpha$ and to make this dependence explicit, we write $c'(0)=T(\alpha)$.  We compute $T(\alpha)$ from \equ{1.2}, differentiating and evaluating in $\tau=0$
\begin{equation}\label{T}
T(\alpha)=\int_0^{+\infty}\frac{\log |1-t^2|}{t^{1+2\alpha}}dt.
\end{equation}
We have to prove that $T$ possesses a unique zero in the interval $(0,1)$.
For this purpose we start proving that
\begin{equation}\label{2.03}
\lim_{\alpha\to1^-}T(\alpha)=-\infty\quad \mbox{and}\quad \lim_{\alpha\to0^+}T(\alpha)=+\infty.
\end{equation}
The first limit follows from the fact that
%
$\log(1-s)\le-s,$ for all $s\in[0,1/4]$, and so
\begin{eqnarray*}
\lim_{\alpha\to 1^-}\int_{0}^{\frac12}\frac{\log (1-t^2)}{t^{1+2\alpha}}dt \leq -\lim_{\alpha\to 1^-} \int_{0}^{\frac12}t^{1-2\alpha}dt=-\infty
\end{eqnarray*}
and the fact that  exists a constant $t_0$ such that
$$
\int_{\frac12}^{+\infty}\frac{\log |1-t^2|}{t^{1+2\alpha}}dt\le t_0,\qquad \mbox{for all } \alpha\in (1/2,1).
$$
The second limit in \equ {2.03} follows from
 \begin{eqnarray*}
 \lim_{\alpha\to0^{+}} \int_{2}^{+\infty}\frac{\log
|1-t^2|}{t^{1+2\alpha}}dt \geq  \log3  \lim_{\alpha\to0^{+}}\int_2^{+\infty}t^{-1-2\alpha}dt= +\infty
\end{eqnarray*}
and the fact that there exists a constant $t_1$ such that
$$
\int_0^{2}\frac{\log
|1-t^2|}{t^{1+2\alpha}}dt\le t_1, \quad \mbox{for all }\alpha\in (0,1/2).
$$
On the other hand we claim that
\begin{equation}\label{2.31}
T'(\alpha)=-2\int_0^{+\infty}\frac{\log|1-t^2|}{t^{1+2\alpha}}\log
tdt<0,\ \ \ \alpha\in(0,1).
\end{equation}
 In fact, since
${\log|1-t^2|}\log t$ is negative  only  for
$t\in(1,\sqrt{2})$, we have
\begin{eqnarray*}
\int_0^{+\infty}\frac{\log|1-t^2|}{t^{1+2\alpha}}\log t dt&>&
\int_0^{\sqrt{2}-1}\frac{\log(1-t^2)}{t^{1+2\alpha}}\log t dt +
\int_1^{\sqrt{2}}\log(t^2-1)\log t dt
\\&\ge&\int_0^{\sqrt{2}-1}\frac{-t^2}{t^{1+2\alpha}}\log t dt +
\int_1^{\sqrt{2}}\log(t-1)\log t dt
\\&=&-\int_0^{\sqrt{2}-1}t^{1-2\alpha}\log t dt +
\int_0^{\sqrt{2}-1}\log(1+t)\log t dt
\\&\ge&-\int_0^{\sqrt{2}-1}t^{1-2\alpha}\log t dt +
\int_0^{\sqrt{2}-1}t\log t dt>0.
\end{eqnarray*}
Then,
 (\ref{2.03}) and (\ref{2.31}) the existence of the desired $\alpha_0\in (0,1)$ with the required properties follows,
 completing  $(i)$ and (\ref{00}) in $(ii)$.

To continue with the proof of our proposition, we study the first limit in (\ref{1.3.2}).
We assume that there exist a
sequence $\alpha_n\in(0,\alpha_0)$ and $\tilde \tau\in(-1,0)$ such
that
$$\lim_{n\to+\infty}\alpha_n=\alpha_0\quad\mbox{  and}\quad \lim_{n\to+\infty}\tau_1(\alpha_n)=\tilde \tau
$$ and so $c(\tilde
\tau)=0$. Moreover $c(0)=0$ and $c'(0)=T(\alpha_0)=0$, contradicting the strict convexity of $c$ given by  (\ref{2.01}).
Next we prove the second limit in (\ref{1.3.2}). We proceed by contradiction, assuming
that there exist a sequence $\{\alpha_n\}\subset(0,1)$ and $
\bar\tau\in(-1,0)$ such that
$$
\lim_{n\to+\infty}\alpha_n=0\quad\mbox{and}\quad \tau_1(\alpha_n)\geq \bar\tau>-1, \quad\mbox{for all } n\in\N.$$
Then there exist $C_1,C_2>0$, depending on $\bar \tau$,  such that
 \begin{eqnarray*}
\int^2_0|\frac{|1-t|^{\tau_1(\alpha_n)}+(1+t)^{\tau_1(\alpha_n)}-2}{t^{1+2\alpha_n}}|dt\le
C_1
\end{eqnarray*}
and
 \begin{eqnarray*}
\lim_{n\to\infty}\int_2^{+\infty}\frac{|1-t|^{\tau_1(\alpha_n)}+(1+t)^{\tau_1(\alpha_n)}-2}{t^{1+2\alpha_n}}dt
\leq -C_2  \lim_{n\to\infty} \int_2^{+\infty}\frac1{t^{1+2\alpha_n}}dt=
-\infty.
\end{eqnarray*}
Then  $c(\tau_1(\alpha_n))\to-\infty$ as $n\to+\infty$, which  is
impossible since $c(\tau_1(\alpha_n))=0.$

We finally prove the last statement of the proposition. Since
$\tau_0(\alpha)\in(-1,0)$ is such that $C(\tau_0(\alpha))=0$ and we have, by definition, that
$$
c(\tau)=C(\tau)+\int_1^{+\infty}\frac{(t-1)^\tau}{t^{1+2\alpha}}dt,
$$
we find that  $c(\tau_0(\alpha))>0$, which together with (\ref{00}), implies
that
$\tau_0(\alpha)\in(-1,\tau_1(\alpha)).$
\hfill$\Box$\\

Next we prove the main proposition in this section, which is on the basis of the construction of super and
sub-solutions.
By hypothesis on the domain $\Omega$ and the manifold $\mathcal{C}$,  there exists
$\delta>0$  such that the distance functions  $d(\cdot)$, to $\partial \Omega$, and $D(\cdot)$, to $\mathcal{C}$, are of class  $C^2$ in $B_\delta $ and $A_\delta$, respectively, and  $dist(A_\delta,
B_\delta)>0$, where $A_\delta=\{x\in \Omega\ |\  D(x)<\delta\}$ and $B_\delta=\{x\in
\Omega\ |\  d(x)<\delta\}$. Now we define the basic function $V_\tau$ as follows
\begin{equation}\label{2.1}
V_\tau(x):=\left\{ \arraycolsep=1pt
\begin{array}{lll}
 D(x)^{\tau},\ & x\in A_\delta\setminus \mathcal{C},\\[2mm]
 d(x)^{2},\ & x\in B_\delta,\\[2mm]
 l(x),\ \ \ \ &
x\in \Omega\setminus  (A_\delta\cup  B_\delta),\\[2mm]
0,\ &x\in\Omega^c,
\end{array}
\right.
\end{equation}
 where  $\tau$ is a parameter in $(-1,0)$ and the function $l$  is  positive  such that
$V_\tau$ is of class $C^2$ in $\R^N\setminus \mathcal{C}$.

\begin{proposition}\label{prop 2.1}
Let $\alpha_0$ and $\tau_1(\alpha)$ be as in Proposition \ref{lm 2.1}.\\
$(i)$\ If $(\alpha,  \tau) \in [\alpha_0,1)\times (-1,0)$
or
$(\alpha,  \tau) \in (0,\alpha_0)\times (-1,\tau_1(\alpha)),$
then there exist $\delta_1\in(0,\delta]$ and $C>1$ such that
$$
\frac1CD(x)^{\tau-2\alpha }\leq-(-\Delta)^{\alpha}V_\tau(x)\leq
CD(x)^{\tau-2\alpha },\ \ x\in A_{\delta_1}\setminus \mathcal{C}.
$$
 $(ii)$\ If
 $(\alpha,  \tau) \in (0,\alpha_0)\times (\tau_1(\alpha),0),$ then there
exist $\delta_1\in(0,\delta]$ and $C>1$ such that
$$
\frac1CD(x)^{\tau-2\alpha }\leq(-\Delta)^{\alpha}V_\tau(x)\leq
CD(x)^{\tau-2\alpha },\ \ x\in A_{\delta_1}\setminus \mathcal{C}.$$
$(iii)$\ If
$(\alpha,  \tau) \in(0,\alpha_0)\times \{\tau_1(\alpha)\},$ then
there exist $\delta_1\in(0,\delta]$ and $C>1$ such that
$$|(-\Delta)^{\alpha}V_\tau(x)|\leq C
D(x)^{\min\{\tau, 2\tau-2\alpha+1\}} , \ \ x\in
A_{\delta_1}\setminus \mathcal{C}.$$
\end{proposition}

This proposition and its proof has many similarities with Proposition 3.2 in \cite{CFQ}, but it has also important differences so we give a complete proof of it.

\medskip

\noindent{\bf Proof.} By compactness of $\mathcal{C}$, we just need to prove that
the corresponding inequality holds in a neighborhood of any point
$\bar x\in\mathcal{C}$ and, without loss of generality, we may  assume $\bar x=0$. For a given $0<\eta\le\delta$, we define
$$Q_\eta=(-\eta,\eta)\times B_\eta\subset\R\times\R^{N-1},$$
where $B_\eta$ denotes the ball centered at the origin and with radius
$\eta$ in $\R^{N-1}$. We observe that $Q_\eta\subset
\Omega.$ Let $\varphi:\R^{N-1}\to \R$ be a $C^2$ function such that
$(z_1,z')\in \mathcal{C}\cap Q_\delta$ if and only if
$z_1=\varphi(z')$. We further assume that $e_1$ is
normal  to $\mathcal{C}$ at $\bar x$ and then there
exists $C>0$ such that $|\varphi(z')|\le C|z'|^2$ for  $|z'|\le
\delta$. Thus, choosing $\eta>0$ smaller if necessary we may assume that
 $|\varphi(z')|<\frac\eta2$ for $|z'|\le \eta$.
 In the proof of our
inequalities, we will consider a generic point along the normal $x=(x_1,0)\in A_{\eta/4}$, with $0<|x_1|< \eta/4$. We observe  that $|x-\bar
x|=D(x)=|x_1|$. By  definition we have
\begin{equation}\label{2.4}
-(-\Delta)^\alpha V_\tau(x)
=\frac12\int_{Q_\eta}\frac{\delta(V_\tau,x,y)}{|y|^{N+2\alpha}}dy+\frac12\int_{\R^N\setminus
Q_\eta}\frac{\delta(V_\tau,x,y)}{|y|^{N+2\alpha}}dy.
\end{equation}
It is not difficult to see that the second integral is bounded by $Cx_1^{\tau}$, for an appropriate constant
$C>0$, so that  we only need to study the
first integral, that from now on we denote by $\frac12E(x_1)$.

Our first goal is to obtain  positive constants $c_1, c_2$  so that lower bound for $E(x_1)$
\begin{equation}\label{2.13}
E(x_1)\ge c_1
c(\tau)|x_1|^{\tau-2\alpha}- c_2|x_1|^{\min\{\tau,2\tau-2\alpha+1\}}
\end{equation}
holds, for all $|x_1|\le \eta/4$.
For this
purpose we assume that $0<\eta\le \delta/2$, then for all $y=(y_1,y')\in
Q_\eta$ we have that $x\pm y\in Q_\delta$, so that
$$
D(x\pm y)\le|x_1\pm y_1-\varphi(\pm y')|,
\quad \mbox{for\ all}\ \ y\in Q_\eta.
$$
From here and the fact that   $\tau\in(-1,0)$, we have that
%
\begin{equation}\label{2.90}
E(x_1)=\int_{Q_\eta}\frac{\delta(V_\tau,x,y)}{|y|^{N+2\alpha}}dy\ge\int_{Q_\eta}\frac{I(y)}{|y|^{N+2\alpha}}dy+\int_{Q_\eta}\frac{J(y)+J(-y)}{|y|^{N+2\alpha}}dy,
\end{equation}
where the functions $I$ and $J$ are defined, for $y\in Q_\eta$, as
\begin{equation}\label{2.7}
I(y)=|x_1-y_1|^\tau+|x_1+y_1|^\tau-2x_1^\tau
\end{equation}
and
\begin{equation}\label{2.8}
J(y)=|x_1+y_1-\varphi(y')|^\tau-|x_1+y_1|^\tau.
\end{equation}
In what follows we assume $x_1>0$ (the case $x_1<0$ is similar).
For the first term of the right hand side in (\ref{2.90}),
 we have
$$
\int_{Q_\eta}\frac{I(y)}{|y|^{N+2\alpha}}dy=x_1^{\tau-2\alpha}\int_{Q_{\frac\eta{x_1}}}\frac{|1-z_1|^\tau+|1+z_1|^\tau-2}{|z|^{N+2\alpha}}dz.
$$
On one hand we have that, for a constant  $c_1$, we have
\begin{eqnarray*}
\int_{\R^N}\frac{|1-z_1|^\tau+|1+z_1|^\tau-2}{|z|^{N+2\alpha}}dz=2c(\tau)\int_{\R^{N-1}}\frac1{(|z'|^2+1)^{\frac{N+2\alpha}2}}dz'
=c_1c(\tau),
\end{eqnarray*}
and, on the other hand, for constants $C_2$ and $C_3$ we have
\begin{eqnarray*}
& &|\int_{-\frac{\eta}{x_1}}^{\frac{\eta}{x_1}}\int_{|z'|\ge \frac{\eta}{x_1}}\frac{|1-z_1|^\tau+|1+z_1|^\tau-2}{|z|^{N+2\alpha}}dz|\\&\le & \int_{-\frac{\eta}{x_1}}^{\frac{\eta}{x_1}}(|1-z_1|^\tau+|1+z_1|^\tau+2)dz_1\int_{|z'|\ge \frac{\eta}{x_1}}\frac{dz'}{|z'|^{N+2\alpha}}\le C_2x_1^{2\alpha}
 \end{eqnarray*}
and
\begin{eqnarray*}
&&|\int_{|z_1|\ge \frac{\eta}{x_1}}\int_{\R^{N-1}}\frac{|1-z_1|^\tau+|1+z_1|^\tau-2}{|z|^{N+2\alpha}}dz|\\
&\le&
2\int_{\frac\eta{x_1}}^{+\infty}\frac{|1-z_1|^\tau+|1+z_1|^\tau+2}{z_1^{1+2\alpha}}dz_1\int_{\R^{N-1}}\frac1{(1+|z'|^2)^{\frac{N+2\alpha}2}}dz'
\le  C_3x_1^{2\alpha}.
\end{eqnarray*}
Consequently, for  an appropriate constant $c_2$
\begin{equation}\label{2.10}
| \int_{Q_\eta}\frac{I(y)}{|y|^{N+2\alpha}}dy -c_1 c(\tau)x_1^{\tau-2\alpha}| \le
c_2 x_1^\tau.
\end{equation}
Next we estimate the second term of the right hand side in (\ref{2.90}). Since
$$\int_{Q_\eta}\frac{J(-y)}{|y|^{N+2\alpha}}dy=\int_{Q_\eta}\frac{J(y)}{|y|^{N+2\alpha}}dy,$$
we only need to estimate
\begin{equation}\label{2.11}
\int_{Q_\eta}\frac{J(y)}{|y|^{N+2\alpha}}dy=\int_{B_\eta}\int^\eta_{-\eta}\frac{|x_1+y_1-\varphi(y')|^\tau
-|x_1+y_1|^\tau}{(y_1^2+|y'|^2)^{\frac{N+2\alpha}2}}dy_1dy'.
\end{equation}
We notice that $|x_1+y_1-\varphi(y')|\ge |x_1+y_1|$ if and only if
$$
\varphi(y')(x_1+y_1-\frac{\varphi(y')}{2})  \le0.
$$
From here and
(\ref{2.11}), we have
\begin{eqnarray*}
\int_{Q_\eta}\frac{J(y)}{|y|^{N+2\alpha}}dy&\ge&\int_{B_\eta}\int^{-x_1+\frac{\varphi_+(y')}2}_{-\eta}\frac{|x_1+y_1-\varphi_+(y')|^\tau
-|x_1+y_1|^\tau}{(y_1^2+|y'|^2)^{\frac{N+2\alpha}2}}dy_1dy'\\
&&+\int_{B_\eta}\int^{\eta}_{-x_1+\frac{\varphi_-(y')}2}\frac{|x_1+y_1-\varphi_-(y')|^\tau
-|x_1+y_1|^\tau}{(y_1^2+|y'|^2)^{\frac{N+2\alpha}2}}dy_1dy'\\
&=&E_1(x_1)+E_2(x_1),
\end{eqnarray*} where
$\varphi_+(y')=\max\{\varphi(y'),0\}$ and
$\varphi_-(y')=\min\{\varphi(y'),0\}$. We only estimate  $E_1(x_1)$ ($E_2(x_1)$ is similar).   Using integration
by parts, we obtain
\begin{eqnarray}& &E_1(x_1) \nonumber \\
&=&\int_{B_\eta}\int^{\frac{\varphi_+(y')}2}_{x_1-\eta}\frac{|y_1-\varphi_+(y')|^\tau
-|y_1|^\tau}{((y_1-x_1)^2+|y'|^2)^{\frac{N+2\alpha}2}}dy_1dy'
\nonumber \\&=&\int_{B_\eta}\int^{0}_{x_1-\eta}\frac{(\varphi_+(y')-y_1)^\tau-(-y_1)^\tau}{((y_1-x_1)^2+|y'|^2)^{\frac{N+2\alpha}2}}dy_1dy'
\nonumber \\
& &+\int_{B_\eta}\int^{\frac{\varphi_+(y')}2}_{0}\frac{(\varphi_+(y')-y_1)^\tau-y_1^\tau}{((y_1-x_1)^2+|y'|^2)^{\frac{N+2\alpha}2}}dy_1dy'
\nonumber \\&=&\frac1{\tau+1}\int_{B_\eta}[\frac{-\varphi_+(y')^{\tau+1}}{(x_1^{2}+|y'|^2)^{\frac{N+2\alpha}2}}
+\frac{(\eta-x_1+\varphi_+(y'))^{\tau+1}-(\eta-x_1)^{\tau+1}}{(\eta^{2}+|y'|^2)^{\frac{N+2\alpha}2}}]dy'
\nonumber \\&&-\frac{N+2\alpha}{\tau+1}\int_{B_\eta}\int^{0}_{x_1-\eta}\frac{(\varphi_+(y')-y_1)^{\tau+1}
-(-y_1)^{\tau+1}}{((y_1-x_1)^2+|y'|^2)^{\frac{N+2\alpha}2+1}}(y_1-x_1)dy_1dy'
\nonumber \\&&+\frac{1}{\tau+1}\int_{B_\eta}[\frac{-2^{-\tau}\varphi_+(y')^{\tau+1}}{((\frac{\varphi_+(y')}{2}-x_1)^{2}+|y'|^2)^{\frac{N+2\alpha}2}}
+\frac{\varphi_+(y')^{\tau+1}}{(x_1^{2}+|y'|^2)^{\frac{N+2\alpha}2}}]dy'
\nonumber \\&&+\frac{N+2\alpha}{\tau+1}\int_{B_\eta}\int^{\frac{\varphi_+(y')}{2}}_{0}\frac{(\varphi_+(y')-y_1)^{\tau+1}
+y_1^{\tau+1}}{((y_1-x_1)^2+|y'|^2)^{\frac{N+2\alpha}2+1}}(y_1-x_1)dy_1dy'
\nonumber \\&\ge&\frac{-2^{-\tau}}{\tau+1}\int_{B_\eta}\frac{\varphi_+(y')^{\tau+1}}{((\frac{\varphi_+(y')}{2}-x_1)^{2}+|y'|^2)^{\frac{N+2\alpha}2}}dy'
\nonumber \\&&+\frac{N+2\alpha}{\tau+1}\int_{B_\eta}\int^{\min\{\frac{\varphi_+(y')}{2},x_1\}}_{0}\frac{(\varphi_+(y')-y_1)^{\tau+1}
+y_1^{\tau+1}}{((y_1-x_1)^2+|y'|^2)^{\frac{N+2\alpha}2+1}}(y_1-x_1)dy_1dy'
\nonumber\\&=&A_1(x_1)+A_2(x_1).\label{cero}
\end{eqnarray}
In order to estimate $A(x_1)$, we split $B_\eta$ in $O=\{y'\in B_\eta: |\frac{\varphi_+(y')}2-x_1|\ge
\frac{x_1}2\}$ and $B_\eta\setminus O$. On one hand we have
\begin{eqnarray*}
\int_{O}\frac{|y'|^{2\tau+2}}{((\frac{\varphi_+(y')}{2}-x_1)^{2}+|y'|^2)^{\frac{N+2\alpha}2}}dy'
 &\le&x_1^{2\tau-2\alpha+1}\int_{B_{\eta/{x_1}}}\frac{|z'|^{2\tau+2}}{(1/4+|z'|^2)^{\frac{N+2\alpha}2}}dz'
 \\&\leq& C(x_1^{2\tau-2\alpha+1}+x_1^\tau).
\end{eqnarray*}
On the other hand,  for $y'\in B_\eta\setminus O$ we have
that $|y'|\ge c_1\sqrt{x_1}$, for some  constant $c_1$, and then
\begin{eqnarray*}\label{dos}
\int_{B_\eta\setminus
O}\frac{|y'|^{2\tau+2}}{((\frac{\varphi_+(y')}{2}-x_1)^{2}+|y'|^2)^{\frac{N+2\alpha}2}}dy'
 &\le&
\int_{B_\eta\setminus
B_{c_1\sqrt{x_1}}}|y'|^{2\tau+2-N-2\alpha}dy'
\nonumber\\&\leq&C(x_1^{\tau-\alpha+\frac12}+1).
\end{eqnarray*}
Thus, for some $C>0$,
\begin{eqnarray}\label{una}
 A_1(x_1)\ge -Cx_1^{\min\{\tau,2\tau-2\alpha+1\}}.
\end{eqnarray}
Next we estimate   $A_2(x_1)$:
\begin{eqnarray*}
A_2(x_1)&\ge&\frac{2(N+2\alpha)}{\tau+1}\int_{B_\eta}\int^{x_1}_{0}
\frac{\varphi_+(y')^{\tau+1}(y_1-x_1)}{((y_1-x_1)^2+|y'|^2)^{\frac{N+2\alpha}2+1}}dy_1dy'
\\&\ge&C\int_{B_\eta}\int^{x_1}_{0}\frac{|y'|^{2\tau+2}(y_1-x_1)}{((y_1-x_1)^2+|y'|^2)^{\frac{N+2\alpha}2+1}}dy_1dy'
\\&\ge&C x_1^{2\tau-2\alpha+1}\int_{B_{\eta/{x_1}}}\int^{1}_{0}\frac{|z'|^{2\tau+2}(z_1-1)}{((z_1-1)^2+|z'|^2)^{\frac{N+2\alpha}2+1}}dz_1dz'
\\&\ge&-C_1 x_1^{\min\{\tau,2\tau-2\alpha+1\}},
\end{eqnarray*}
for some $C,C_1>0$. From here, \equ{cero} and  \equ{una} we obtain, for some $C>0$
$$
E_1(x_1)\ge
-Cx_1^{\min\{\tau,2\tau-2\alpha+1\}}.$$
Using the similar estimate for $E_2(x_1)$, we obtain
\begin{equation}\label{2.00}
\int_{Q_\eta}\frac{J(y)+J(-y)}{|y|^{N+2\alpha}}dy\ge-Cx_1^{\min\{\tau,2\tau-2\alpha+1\}}.
\end{equation}
Thus, from  (\ref{2.90}), (\ref{2.10}), (\ref{2.00}) and noticing that these inequalities also hold with $x_1<0$ with the obvious changes, we conclude the lower bound for $E(x_1)$ we gave in \equ{2.13}.
Our second goal is to get an upper bound for $E(x_1)$ and for this, we first recall Lemma 3.1 in
\cite{CFQ} to obtain
$$
D(x\pm y)^\tau\le(x_1\pm y_1-\varphi(y'))^\tau(1+C|y'|^2),\,\,\mbox{for all}\quad  |x_1|\le \eta/4, y=(y_1,y')\in
Q_\eta.
$$
From here we see that
\begin{eqnarray}\label{2.9}
E(x_1)&\le&\int_{Q_\eta}\frac{I(y)}{|y|^{N+2\alpha}}dy+\int_{Q_\eta}\frac{J(y)+J(-y)}{|y|^{N+2\alpha}}dy\nonumber
\\ &&+C\int_{Q_\eta}\frac{I(y)+J(y)+J(-y)}{|y|^{N+2\alpha}}|y'|^2dy.
\end{eqnarray}
We denote by
$E_3(x_1)$ the third integral above. The first integral was studied in  \equ{2.10}, so we study the second integral and that we only need to consider the term $J(y)$, since the other is completely analogous.
We see that
 $|x_1+y_1-\varphi(y')|\le|x_1+y_1|$ if and only if
$$
\varphi(y')(x_1+ y_1-\frac{\varphi(y')}2)\ge 0.
$$
As before, we will consider only the case $x_1>0$, since the other one is analogous. From (\ref{2.11}) we have
 \begin{eqnarray*}
\int_{Q_\eta}\frac{J(y)}{|y|^{N+2\alpha}}dy&\le&\int_{B_\eta}\int^{-x_1+\frac{\varphi_-(y')}2}_{-\eta}\frac{|x_1+y_1-\varphi_-(y')|^\tau
-|x_1+y_1|^\tau}{(y_1^2+|y'|^2)^{\frac{N+2\alpha}2}}dy_1dy'
\\&&+\int_{B_\eta}\int^{\eta}_{-x_1+\frac{\varphi_+(y')}2}\frac{|x_1+y_1-\varphi_+(y')|^\tau
-|x_1+y_1|^\tau}{(y_1^2+|y'|^2)^{\frac{N+2\alpha}2}}dy_1dy'
\\&=&F_1(x_1)+F_2(x_1).
\end{eqnarray*}
 Next we estimate
$F_1(x_1)$ ($F_2(x_1)$ is similar), using integration by parts
\begin{eqnarray*}
& & F_1(x_1)\\
&=&\int_{B_\eta}\int_{x_1-\eta}^{\frac{\varphi_-(y')}2}\frac{|y_1-\varphi_-(y')|^\tau
-|y_1|^\tau}{((y_1-x_1)^2+|y'|^2)^{\frac{N+2\alpha}2}}dy_1dy'
\\&=&\int_{B_\eta}\int^{\varphi_-(y')}_{x_1-\eta}\frac{(\varphi_-(y')-y_1)^\tau-(-y_1)^\tau}{((y_1-x_1)^2+|y'|^2)^{\frac{N+2\alpha}2}}dy_1dy'
\\&&+\int_{B_\eta}\int^{\frac{\varphi_-(y')}2}_{\varphi_-(y')}\frac{(y_1-\varphi_-(y'))^\tau-(-y_1)^\tau}{((y_1-x_1)^2+|y'|^2)^{\frac{N+2\alpha}2}}dy_1dy'
\\&=&\frac1{\tau+1}\int_{B_\eta}[\frac{(-\varphi_-(y'))^{\tau+1}}{((x_1-\varphi_-(y'))^{2}+|y'|^2)^{\frac{N+2\alpha}2}}
+\frac{(\eta-x_1+\varphi_-(y'))^{\tau+1}-(\eta-x_1)^{\tau+1}}{(\eta^{2}+|y'|^2)^{\frac{N+2\alpha}2}}]dy'
\\&&-\frac{N+2\alpha}{\tau+1}\int_{B_\eta}\int^{\varphi_-(y')}_{x_1-\eta}\frac{(\varphi_-(y')-y_1)^{\tau+1}
-(-y_1)^{\tau+1}}{((y_1-x_1)^2+|y'|^2)^{\frac{N+2\alpha}2+1}}(y_1-x_1)dy_1dy'
\\&&+\frac{1}{\tau+1}\int_{B_\eta}[\frac{2^{-\tau}(-\varphi_-(y'))^{\tau+1}}{((\frac{\varphi_-(y')}{2}-x_1)^{2}+|y'|^2)^{\frac{N+2\alpha}2}}
+\frac{-(-\varphi_-(y'))^{\tau+1}}{((x_1-\varphi_-(y'))^{2}+|y'|^2)^{\frac{N+2\alpha}2}}]dy'
\\&&+\frac{N+2\alpha}{\tau+1}\int_{B_\eta}\int^{\frac{\varphi_-(y')}{2}}_{\varphi_-(y')}\frac{(y_1-\varphi_-(y'))^{\tau+1}
+(-y_1)^{\tau+1}}{((y_1-x_1)^2+|y'|^2)^{\frac{N+2\alpha}2+1}}(y_1-x_1)dy_1dy'
\\&\le&
\frac{1}{\tau+1}\int_{B_\eta}\frac{2^{-\tau}(-\varphi_-(y'))^{\tau+1}}{((\frac{\varphi_-(y')}{2}-x_1)^{2}+|y'|^2)^{\frac{N+2\alpha}2}}dy'=B(x_1).
\end{eqnarray*}
Since
$(\frac{\varphi_-(y')}{2}-x_1)^{2}\ge x_1^2$, we have
\begin{eqnarray*}
B(x_1)&\le&
\frac{2^{-\tau}}{\tau+1}\int_{B_\eta}\frac{(-\varphi_-(y'))^{\tau+1}}{(x_1^{2}+|y'|^2)^{\frac{N+2\alpha}2}}dy'
\\&\le& C\int_{B_\eta}\frac{|y'|^{2\tau+2}}{(x_1^{2}+|y'|^2)^{\frac{N+2\alpha}2}}dy'
\le Cx_1^{\min\{\tau,2\tau-2\alpha+1\}},
\end{eqnarray*}
for some $C>0$ independent of $x_1$. Thus we have obtained that
\begin{equation}\label{tres}
F_1(x_1)\le Cx_1^{\min\{\tau,2\tau-2\alpha+1\}}.
\end{equation}
Similarly, we can get an analogous estimate for $F_2(x_1)$ and these two estimates imply
\begin{equation}\label{2.00i}
\int_{Q_\eta}\frac{J(y)+J(-y)}{|y|^{N+2\alpha}}dy\le Cx_1^{\min\{\tau,2\tau-2\alpha+1\}}.
\end{equation}
Finally we obtain
\begin{eqnarray*}
\int_{Q_\eta}\frac{I(y)}{|y|^{N+2\alpha}}|y'|^2dy&=&
x_1^{\tau-2\alpha+2}\int_{Q_{\frac\eta{x_1}}}\frac{|1-z_1|^\tau+|1+z_1|^\tau-2}{|z|^{N+2\alpha}}|z'|^2dz
\\&\le&Cx_1^{\min\{\tau,\tau-2\alpha+2\}}
\end{eqnarray*}
and, in a similar way,
$$\int_{Q_\eta}\frac{J(y)|y'|^2}{|y|^{N+2\alpha}}dy\le Cx_1^{\min\{\tau, 2\tau-2\alpha+1\}}.
$$
From the last two inequalities we obtain
\begin{equation}\label{3.0.4}
E_3(x_1)\le Cx_1^{\min\{\tau, 2\tau-2\alpha+1\}}.
\end{equation}
Then, taking into account  (\ref{2.9}),  (\ref{2.10}), \equ{2.00i}, (\ref{3.0.4})
 and considering also the case $x_1<0$, we obtain
\begin{equation}\label{2.13i}
E(x_1)\le c_1
c(\tau)|x_1|^{\tau-2\alpha}+ c_2|x_1|^{\min\{\tau,2\tau-2\alpha+1\}}.
\end{equation}
From inequalities \equ{2.13},  \equ{2.13i}  and
Proposition \ref{lm 2.1}
the result follows. \hfill$\Box$

\setcounter{equation}{0}
\section{Existence of large solution}
This section is devoted to use Proposition \ref{prop 2.1} to prove
the existence of solution of problem (\ref{inside}). To this
purpose, our main goal is to construct appropriate  sub-solution and
super-solution of problem (\ref{inside}) under the hypotheses of
Theorem \ref{th 1.2} $(i)$, $(ii)$ and Theorem \ref{th 1.1} $(i)$.

We begin with a simple lemma that reduces the problem to find them only in $A_\delta \setminus \mathcal{C}$.

\begin{lemma}\label{lemall}
Let $U$ and $W$  be classical ordered super and sub-solution  of \equ{inside}  in the sub-domain  $A_\delta\setminus \mathcal{C}$.
Then there exists $\lambda$ large such that
$U_\lambda=U+\lambda\bar V$ and  $ W_\lambda=W-\lambda\bar V$, are ordered super and sub-solution  of \equ{inside},
where $\bar V$ is the solution of
\begin{equation} \label{3.1}
\left\{ \arraycolsep=1pt
\begin{array}{lll}
 (-\Delta)^{\alpha} \bar V(x)=1,\ \ \ \ &
x\in\Omega,\\[2mm]
\bar V(x)=0,\ &x\in\Omega^c.
\end{array}
\right.
\end{equation}
\end{lemma}
\begin{remark}Here $U,W:\RR^N\to \R$ are classical ordered of  super and sub-solution of \equ{inside}  in the sub-domain  $A_\delta\setminus \mathcal{C}$
if $U$ satisfies
$$ (-\Delta)^\alpha U+|U|^{p-1}U\geq 0 \quad \mbox{in}\quad A_\delta\setminus \mathcal{C}$$
and $W$ satisfies the reverse inequality. Moreover, they satisfy
$$
U\geq W\ \ {\rm{in}}\ \Omega\setminus\mathcal{C},\quad  \liminf_{x\in\Omega\setminus\mathcal{C},x\to\mathcal{C}} W(x)=+\infty,\quad U=W=0\ \ {\rm{in}}\ \Omega^c.$$

\end{remark}
\noindent {\bf Proof.}
Notice that by the maximum principle $\bar V$ is nonnegative in $\Omega$, therefore $U_\lambda\geq U$ and $W_\lambda\leq W$, so they are still ordered.
In addition $U_\lambda$ satisfies $$(-\Delta)^\alpha U_\lambda+|U_\lambda|^{p-1}U_\lambda\ge (-\Delta)^\alpha U+|U|^{p-1}U+\lambda>0,\quad\mbox{in}\quad\Omega\setminus \mathcal{C}.$$
This inequality holds because of our  assumption in $A_\delta\setminus \mathcal{C}$ and the fact that $ (-\Delta)^\alpha U+|U|^{p-1}U$ is continuous in $\Omega\setminus{A_\delta}$ and
by  taking $\lambda$ large enough.

By the same type of arguments we find that $W_\lambda$ is a sub-solution.  \hfill$\Box$

\medskip
\noindent {\bf Proof of existence results in Theorem \ref{th 1.2}
$(i)$ and Theorem \ref{th 1.1} $(i)$.}
We define
\begin{equation}\label{3.4.1}
U_\lambd(x)=\lambd V_\tau(x) \ {\rm{ and}}\ \
W_{\lambd}(x)=\lambd V_\tau(x),\
x\in\R^N\setminus\mathcal{C},
\end{equation}
where $V_\tau$ is defined in (\ref{2.1}) with
$\tau=-\frac{2\alpha}{q-1}$

\textbf{1. $U_\lambd$ is Super-solution.} By hypotheses of Theorem \ref{th 1.2} $(i)$ and Theorem \ref{th 1.1}
$(i)$, we notice that  $$\tau\in(-1,0),\quad\rm{for}\
\alpha\in[\alpha_0,1),$$ $$ \tau\in
(-1,\tau_1(\alpha)),\quad\rm{for}\
 \alpha\in(0,\alpha_0)$$ and $\tau p=\tau-2\alpha $, then we use
Proposition \ref{prop 2.1} part $(i)$ to obtain that there exist
$\delta_1\in(0,\delta]$ and $C>1$ such that
\begin{eqnarray*}
(-\Delta)^\alpha U_\lambd(x)+U^p_\lambd(x)\geq -C\lambd
D(x)^{\tau-2\alpha }+\lambd^p D(x)^{\tau p},\quad x\in
A_{\delta_1}\setminus\mathcal{C}.
\end{eqnarray*}
Then there exist $\lambd_1>1$ such that for $\lambd\geq
\lambd_1$, we have
$$(-\Delta)^\alpha U_\lambd(x)+U^p_\lambd(x)\geq 0,\ x\in
A_{\delta_1}\setminus\mathcal{C}.$$

\noindent
{\bf 2. $W_{\lambd}$ is Sub-solution. }  We use  Proposition \ref{prop 2.1} part
$(i)$ to obtain that there exist $\delta_1\in(0,\delta]$ and $C>1$
such that for $x\in A_{\delta_1}\setminus\mathcal{C}$, we have
\begin{eqnarray*}
(-\Delta)^\alpha
 W_{\lambd}(x)+ |W_{\lambd}|^{p-1}W_{\lambd}(x)&\leq&
-\frac\lambd C D(x)^{\tau-2\alpha}+\lambd^p D(x)^{\tau p}\\&\leq&
(-\frac\lambd C+\lambd^p)D(x)^{\tau-2\alpha}.
\end{eqnarray*}
Then there exists $\lambd_3\in (0,1)$ such
that for all $\lambd\in(0,\lambd_3)$, it has
$$(-\Delta)^\alpha  W_{\lambd}(x)+|W_{\lambd}|^{p-1}W_{\lambd}(x)\leq 0,\ x\in
A_{\delta_1}\setminus\mathcal{C}.$$
To conclude the proof we use  Lemma \ref{lemall}  and Proposition \ref{prop 2.1}.  \hfill $\Box$

\medskip

\noindent{\bf Proof of Theorem \ref{th 1.2} $(ii)$.}  For any given
$t>0$, we denote
$$
U(x)=t V_{\tau_1(\alpha)}(x),\quad
x\in\R^N\setminus\mathcal{C},$$
$$ W_{\mu}(x)=t V_{\tau_1(\alpha)}(x)-\mu V_{\bar \tau}(x),\quad x\in\R^N\setminus\mathcal{C}$$ where $\bar
\tau=\min\{\tau_1(\alpha)p+2\alpha,\frac12\tau_1(\alpha)\}<0$. By (\ref{observartion}), we
have
\begin{equation}\label{tau bar}
\bar \tau\in(\tau_1(\alpha),0),\ \bar \tau-2\alpha<
\min\{\tau_1(\alpha),2\tau_1(\alpha)-2\alpha +1\}\ \rm{and}\  \bar
\tau-2\alpha<\tau_1(\alpha)p.
\end{equation}

\noindent\textbf{1. $U$ is Super-solution.}  We use Proposition
\ref{prop 2.1} $(iii)$ to obtain that for any  $x\in
A_{\delta_1}\setminus\mathcal{C}$,
\begin{eqnarray*}
(-\Delta)^\alpha U(x)+U^p(x)\geq -Ct
 D(x)^{\min\{\tau_1(\alpha),2\tau_1(\alpha)-2\alpha +1\}
}+t^p D(x)^{\tau_1(\alpha) p},
\end{eqnarray*}
together with
$\tau_1(\alpha)p<\min\{\tau_1(\alpha),2\tau_1(\alpha)-2\alpha +1\}$,
then  there exists $\delta_2\in(0,\delta_1]$ such that
$$(-\Delta)^\alpha U(x)+U^p(x)\geq 0,\quad x\in
A_{\delta_2}\setminus\mathcal{C}.$$

\noindent
\textbf{2. $W_{\mu}$ is Sub-solution.}  We use Proposition \ref{prop 2.1} $(ii)$ and
$(iii)$ to obtain that  for $x\in A_{\delta_1}\setminus\mathcal{C}$,
\begin{eqnarray*}(-\Delta)^\alpha
 W_{\mu}(x)+ |W_{\mu}|^{p-1}W_{\mu}(x)&\leq&
Ct
 D(x)^{\min\{\tau_1(\alpha),2\tau_1(\alpha)-2\alpha+1\}
}\\&&-\frac\mu C D(x)^{\bar \tau-2\alpha }+t^p D(x)^{\tau_1(\alpha)
p}.
\end{eqnarray*}
Then there exists $\delta_2\in(0,\delta_1]$ such that for any $\mu\ge1$,
we have
$$(-\Delta)^\alpha  W_{\mu}(x)+|W_{\mu}|^{p-1}W_{\mu}(x)\leq 0,\ x\in
A_{\delta_2}\setminus\mathcal{C}.$$
To conclude the proof we use  Lemma \ref{lemall}  and Proposition \ref{prop 2.1}. \hfill $\Box$

\section{Uniqueness and nonexistence}
 We prove the uniqueness statement by contradiction. Assume that $u$ and $v$ are solutions of problem (\ref{inside})
 satisfying (\ref{0.1.2}). Then there exist $C_0\ge1$ and $\bar
 \delta\in(0,\delta)$ such that
\begin{equation}\label{4.1.3}
 \frac 1{C_0}\leq v(x)D(x)^{-\tau},\ u(x)D(x)^{-\tau} \leq C_0,\ \ \forall x\in
  A_{\bar\delta}\setminus\mathcal{C},
\end{equation}
where $\tau=-\frac{2\alpha}{p-1}$. We denote
\begin{equation}\label{4.1.4}
\mathcal{A}=\{x\in  \Omega\setminus\mathcal{C}\ |\
u(x)>v(x)\}.
\end{equation}
 Then
$\mathcal{A}$ is open and $\mathcal{A}\subset \Omega$. Then the
uniqueness in Theorem \ref{th 1.1} $(i)$ and Theorem \ref{th 1.2}
$(i)$ is a consequence of the following result:
\begin{proposition}\label{th 3.1}
Under the hypotheses of Theorem \ref{th 1.1} $(i)$ and Theorem
\ref{th 1.2} $(i)$, we have
$$\mathcal{A}=\O.$$
\end{proposition}
{\bf Proof.} The procedure of proof is similar as Section\S 5 in
\cite{CFQ}, noting that we need to replace $d(x)$ by $D(x)$ and $\partial \Omega$ by $\mathcal{C}$
.\hfill$\Box$\\

From Proposition \ref{th 3.1}, we can prove  uniqueness part  in Theorem \ref{th 1.2} $(i)$ and
Theorem \ref{th 1.1} $(i)$ .

The final goal in this note is to consider the nonexistence of
 solutions of problem (\ref{inside}) under the hypotheses of Theorem \ref{th 1.2} $(iii)$ and Theorem \ref{th 1.1} $(ii)$.

\begin{proposition}\label{th 4}
Under the hypotheses of Theorem \ref{th 1.2} $(iii)$ and Theorem
\ref{th 1.1} $(ii)$, we assume that $U_1$ and $U_2$ are both
sub-solutions (or both super-solutions) of (\ref{inside})  satisfying
that $U_1=U_2=0$ in $\Omega^c$ and
\begin{eqnarray*}
0&<&\liminf_{x\in\Omega\setminus\mathcal{C},\
x\to\mathcal{C}}U_1(x)D(x)^{-\tau}\le\limsup_{x\in\Omega\setminus\mathcal{C},\
x\to\mathcal{C}}U_1(x)D(x)^{-\tau}\\&<&\liminf_{x\in\Omega\setminus\mathcal{C},\
x\to\mathcal{C}}U_2(x)D(x)^{-\tau}\le\limsup_{x\in\Omega\setminus\mathcal{C},\
x\to\mathcal{C}}U_2(x)D(x)^{-\tau}<+\infty,
\end{eqnarray*}
for  $\tau\in(-1,0)$. For the case $\tau p>\tau-2\alpha$, we further
assume that\\
$(i)$ if $U_1,U_2$ are  sub-solutions, there exist $C>0$ and $\tilde
\delta>0$,
\begin{equation}\label{5.1}
(-\Delta)^\alpha U_2(x)\le -CD(x)^{\tau-2\alpha},\quad x\in
A_{\tilde\delta}\setminus\mathcal{C};
\end{equation}
or\\ $(ii)$ if $U_1,U_2$ are  super-solutions, there exist $C>0$ and
$\tilde \delta>0$,
\begin{equation}\label{5.2}
(-\Delta)^\alpha U_1(x)\ge CD(x)^{\tau-2\alpha},\quad x\in
A_{\tilde\delta}\setminus\mathcal{C}.
\end{equation}
Then there doesn't exist any  solution $u$ of (\ref{inside}) such
that
\begin{equation}\label{4.1}
\limsup_{x\in\Omega\setminus\mathcal{C},\
x\to\mathcal{C}}\frac{U_1(x)}{u(x)}<1<\liminf_{x\in\Omega\setminus\mathcal{C},\
x\to\mathcal{C}}\frac{U_2(x)}{u(x)}.
\end{equation}
\end{proposition}
{\bf Proof.} The proof is similar as Proposition 6.1 in
\cite{CFQ},  noting again that  we need to replace $d(x)$ by $D(x)$ and $\partial \Omega$ by $\mathcal{C}$
.\hfill$\Box$\\[1mm]

With the help of Proposition \ref{prop 2.1},  for given $t_1>t_2>0$,
we construct two sub-solutions (or both super-solutions) $U_1$ and
$U_2$ of (\ref{inside}) such that
$$\lim_{x\in\Omega\setminus\mathcal{C},x\to\mathcal{C}}U_1(x)D(x)^{-\tau}=t_1,
\
\lim_{x\in\Omega\setminus\mathcal{C},x\to\mathcal{C}}U_2(x)D(x)^{-\tau}=t_2.$$
So what we have to do is to prove that for any $t>0$, we can
construct super-solution (sub-solution) of problem
(\ref{inside}).

\noindent{\bf Proof of  Theorem \ref{th 1.2} $(iii)$ and Theorem
\ref{th 1.1} $(ii)$.} We divide our proof of the nonexistence
results into several cases under the assumption $p>1$.

\noindent \textbf{Zone 1:} We consider $\tau\in(\tau_1(\alpha),0)$
and $ \alpha\in (0,\alpha_0). $  By Proposition \ref{prop 2.1}
$(ii)$, there exists
 $\delta_1>0$ such that
\begin{equation}\label{4.4}
(-\Delta)^{\alpha}V_\tau(x)\ge\frac1CD(x)^{\tau-2\alpha},\ \ x\in
A_{\delta_1}\setminus\mathcal{C}.
\end{equation}
Since $V_\tau$ is $C^2$  in $\Omega\setminus\mathcal{C}$, then there
exists $C>0$ such that
\begin{equation}\label{4.5}
|(-\Delta)^\alpha V_\tau(x)|\leq C,\ \ x\in\Omega\setminus
A_{\delta_1}.
\end{equation}
Let $\bar U:=V_\tau+C \bar V\quad\rm{in}\ \
\R^N\setminus\mathcal{C}$, then we have $\bar U>0$ in
$\Omega\setminus\mathcal{C}$,
$$
(-\Delta)^\alpha \bar U\ge 0\ \ {\rm{in}}\ \ \Omega\setminus\mathcal{C}
\quad{\rm{and}}\ \ (-\Delta)^\alpha \bar U(x)\ge
\frac1CD(x)^{\tau-2\alpha},\ \ x\in
A_{\delta_1}\setminus\mathcal{C}.
$$
Then, we have that $t\bar U$ is
super-solution of (\ref{inside}) for any $t>0$. Using Proposition
\ref{th 4}, we see that there is no solution of
(\ref{inside}) satisfying (\ref{0.1.3}).\\[1mm]
\noindent \textbf{Zone 2:} We consider $\tau-2\alpha<\tau p$ and
$$\tau\in\left\{ \arraycolsep=1pt
\begin{array}{lll}
 (-1,0),\ & \alpha\in [\alpha_0,1),\\[2mm]
(-1,\tau_1(\alpha)),\quad & \alpha\in (0,\alpha_0).
\end{array}
\right.
$$
 Let us define
$$W_{\mu,t}=tV_\tau-\mu \bar V\quad{\rm{in}}\ \ \R^N\setminus\mathcal{C},$$
where $t,\mu>0$. By Proposition \ref{prop 2.1} $(i)$, for $x\in
A_{\delta_1}\setminus\mathcal{C}$,
\begin{eqnarray*}
 (-\Delta)^\alpha
 W_{\mu,t}(x)+
 |W_{\mu,t}|^{p-1}W_{\mu,t}(x)\leq-\frac tC D(x)^{\tau-2\alpha}+t^pD(x)^{\tau
 p}.
\end{eqnarray*}
 For any fixed $t>0$, there exists $\delta_2\in(0,\delta_1]$, for all
 $\mu\ge0$,
\begin{equation}\label{4.7}(-\Delta)^\alpha
 W_{\mu,t}(x)+
 |W_{\mu,t}|^{p-1}W_{\mu,t}(x)\leq0,\quad A_{\delta_2}\setminus\mathcal{C}.\end{equation}
To consider $x\in \Omega\setminus A_{\delta_2}$, in fact, there
exists $C_1>0$ such that
$$t|(-\Delta)^\alpha V_\tau(x)|+t^pV_\tau^p(x)\le C_1,\quad x\in\Omega\setminus
A_{\delta_2}$$ and
\begin{eqnarray*}
 (-\Delta)^\alpha
 W_{\mu,t}(x)+ |W_{\mu,t}|^{p-1}W_{\mu,t}(x)\leq C_1t-\mu,\quad x\in\Omega\setminus
A_{\delta_2}
\end{eqnarray*}
For given $t>0$, there exists $\mu(t)>0$  such that
\begin{equation}\label{4.9}
 (-\Delta)^\alpha
 W_{\mu(t),t}(x)+|W_{\mu,t}|^{p-1}W_{\mu(t),t}(x)\leq0,\ \ x\in\Omega\setminus A_{\delta_2}.
\end{equation}
Therefore, together with (\ref{4.7}) and (\ref{4.9}), for any given
$t>0$, there sub-solutions $W_{\mu(t),t}$ of problem (\ref{inside})
and by Proposition \ref{th 4}, we see that there is no solution $u$ of
(\ref{inside}) satisfying (\ref{0.1.3}).\\[1mm]
\noindent\textbf{Zone 3:} We consider $\tau-2\alpha>\tau p$ and
$$\tau\in\left\{ \arraycolsep=1pt
\begin{array}{lll}
 (-1,0),\ & \alpha\in [\alpha_0,1),\\[2mm]
(-1,\tau_1(\alpha)),\quad & \alpha\in (0,\alpha_0).
\end{array}
\right.
$$

We denote that
$$U_{\mu,t}=tV_\tau+\mu \bar V\quad\rm{in}\ \ \R^N\setminus\mathcal{C},$$
where $t,\mu>0$. Here $U_{\mu,t}>0$ in $\Omega\setminus\mathcal{C}$.
By Proposition \ref{prop 2.1} $(i)$,
\begin{eqnarray*}
 (-\Delta)^\alpha
 U_{\mu,t}(x)+
 U^p_{\mu,t}(x)\geq-Ct D(x)^{\tau-2\alpha}+t^pD(x)^{\tau p},\quad x\in A_{\delta_1}\setminus\mathcal{C}.
\end{eqnarray*}
For any fixed $t>0$, there exists $\delta_2\in(0,\delta_1]$, for all
 $\mu\ge0$,
\begin{equation}\label{4.701}
(-\Delta)^\alpha  U_{\mu,t}(x)+ U^p_{\mu,t}(x)\geq0,\quad x\in
A_{\delta_2}\setminus\mathcal{C}.
\end{equation}
For $x\in \Omega\setminus A_{\delta_2}$, we see that
$(-\Delta)^\alpha V_\tau$ is bounded  and
\begin{eqnarray*}
 (-\Delta)^\alpha
 U_{\mu,t}(x)+
 U^p_{\mu,t}(x)\geq-Ct+\mu.
\end{eqnarray*}
For given $t>0$, there exists $\mu(t)>0$  such that
\begin{equation}\label{4.801}
 (-\Delta)^\alpha
 U_{\mu(t),t}(x)+
 U^p_{\mu(t),t}(x)\geq0,\ \ x\in\Omega\setminus A_{\delta_2}.
\end{equation}
Combining with (\ref{4.701}) and (\ref{4.801}), we have that for any
$t>0$, there exists $\mu(t)>0$ such that $$ (-\Delta)^\alpha
 U_{\mu(t),t}(x)+
 U^p_{\mu(t),t}(x)\geq0,\ \ \ x\in\Omega\setminus\mathcal{C}.$$
Therefore, for any given $t>0$, there is a super-solution $U_{\mu(t),t}$
of problem (\ref{inside}) and by Proposition \ref{th 4}, we see that there is no solution of (\ref{inside}) satisfying (\ref{0.1.3}).

We see that Zones 1 and 2 cover Theorem \ref{th 1.2} part $(iii)$ a) since $\tau>-2\alpha/(p-1).$
 From Zones 1, 2 and 3 we cover Theorem \ref{th 1.2} part $(iii)$ b) since $\tau_1(\alpha)>2\alpha/(p-1)$.
Moreover, from Zone 1 to Zone 3, we cover the parameters in part $(iii)$ c) of Theorem \ref{th 1.2}, since $\tau_1(\alpha)<2\alpha/(p-1)$.
Finally Theorem \ref{th 1.1} part ii) can be obtained  from Zone 2 and Zone 3.
This complete the proof.
 \hfill$\Box$

\medskip

\noindent
{\bf Acknowledgements}. The authors thanks Peter Bates for proposing the problem.
H.C. was partially supported by Conicyt Ph.D. scholarship. P.F. was partially supported by Fondecyt \# 1110291 and Programa BASAL-CMM U. de Chile.
A.Q. was partially supported by Fondecyt \# 1110210 and Programa BASAL-CMM U. de Chile.

\end{document}